\documentclass[a4paper, 12pt]{extarticle}
\usepackage{geometry}
\usepackage{graphics}
\usepackage{dcolumn}
\usepackage{bm}
\usepackage{amsmath}
\usepackage{hyperref}
\hypersetup{colorlinks=true, urlcolor=blue}
\usepackage{color}
\usepackage{tikz}
\usepackage{multirow}
\usepackage{cite}
\usepackage{float}
\usepackage[caption = false]{subfig}
\usepackage{newlfont}
\usepackage{amssymb}
\usepackage{amsfonts}
\usepackage{amsthm}
\usepackage{graphicx}
\usepackage{graphics,wrapfig,times}
\usepackage{bm}

\def \be{\begin{equation}}
\def \ee{ \end{equation} }
\def \bea{\begin{eqnarray}}
\def \eea{\end{eqnarray}}

\begin{document}

\definecolor{red}{rgb}{1,0,0}
\title{Relations between \(e\), \(\pi\), golden ratios and \(\sqrt{2}\)}
\author{Asutosh Kumar\\
\small P. G. Department of Physics, Gaya College, Magadh University, Rampur, Gaya 823001, India\\ 
\small Vaidic and Modern Physics Research Centre, Bhagal Bhim, Bhinmal, Jalore 343029, India\\ \small (asutoshk.phys@gmail.com)}
\date{}

\maketitle

\textbf{Abstract.} 
We write out relations between the base of natural logarithms ($e$), the ratio of the circumference of a circle to its diameter ($\pi$), the golden ratios ($\Phi_p$) of the additive \(p\)-sequences, and the ratio of the diagonal of a square to its side ($\sqrt{2}$).
An additive \(p\)-sequence is a natural extension of the Fibonacci sequence in which every term is the sum of {\it \(p\) previous} terms given \(p \ge 1\) initial values called {\it seeds}.

\section{Introduction}
Euler's identity (or Euler's equation) is given as
\begin{equation}\label{eq:euler-iden}
e^{i\pi }+1=0,
\end{equation}
where
$e = 2.718 \cdots$ is the base of natural logarithms,
$i := \sqrt{-1}$ is the imaginary unit of complex numbers, and
$\pi = 3.1415 \cdots$ is the ratio of the circumference of a circle to its diameter. 
It is a special case of Euler's formula, $e^{i \theta}=\cos \theta + i \sin \theta$, for $\theta = \pi$. 
This expresses a deep mathematical beauty \cite{nahin, stipp, wilson} as it involves three of the basic arithmetic operations: addition/subtraction, multiplication/division, and exponentiation/logarithm, and five fundamental mathematical constants: $0$ (the additive identity), $1$ (the multiplicative identity), $e$ (Euler's number), $i$ (the imaginary unit), and $\pi$ (the fundamental circle constant).

As Euler's identity is an example of mathematical elegance, further generalizations of similar-type have been discovered. 
\begin{itemize}
\item The $n^{th}$ roots of unity ($n > 1$) add up to zero.
\begin{equation}\label{eq:roots-unity}
\sum^{n-1}_{k=0} e^{2i\pi \frac{k}{n}}=0.
\end{equation}
It yields Euler's identity (\ref{eq:euler-iden}) when $n=2$.

\item For quaternions \cite{jia}, with the basis elements $\{i, j, k\}$ and real numbers $a_n$ such that $a_1^2 + a_2^2 + a_3^2 = 1$,
\begin{equation}\label{eq:octo}
e^{\left(a_1i +a_2j + a_3k \right) \pi} + 1=0.
\end{equation}

\item For octonions, with the basis elements $\{i_1, i_2, \cdots, i_7\}$ and real numbers $a_n$ such that $a_1^2 + a_2^2 + \cdots + a_7^2 = 1$,
\begin{equation}\label{eq:octo}
e^{\left(\sum^{7}_{k=1} a_ki_k \right) \pi} + 1=0.
\end{equation}
\end{itemize}

In this article, motivated by Euler's identity and its generalizations, we give relations between the base of natural logarithms ($e$), the ratio of the circumference of a circle to its diameter ($\pi$), and the golden ratios ($\Phi_p$) of the additive \(p\)-sequences. 

\section{Additive \(p\)-sequences}
An additive \(p\)-sequence \cite{asutosh-gr} is a natural extension of the Fibonacci sequence \cite{koshy2001} in which every term is the sum of {\it \(p\) previous} terms given \(p \ge 1\) initial values called {\it seeds} $(s_0,s_1,\cdots,s_{p-1})$ such that $t_0=s_0,~t_1=s_1,\cdots,t_{p-1}=s_{p-1}$, and
\be t_n(p) := t_{n-1}(p) + t_{n-2}(p) + \cdots + t_{n-p}(p) = \sum_{k=n-p}^{n-1} t_k(p). \ee 
This can be equivalently rewritten as
\be t_{n+p}(p) := t_{n+p-1}(p) + t_{n+p-2}(p) + \cdots + t_{n}(p). \ee 
Varying the values of seeds, it is possible to construct an infinite number of \(p\)-sequences.

By definition of $t_n(p)$, we have
$t_{n+1}(p) >  t_n(p)$ and $t_{n+1}(p) = 2t_n(p) -t_{n-p}(p) < 2t_n(p)$ \cite{asutosh-gr}. Hence
\be 1 < \Phi_p < 2. \ee

The limiting ratio value $\left( \lim_{n \rightarrow \infty} \frac{t_{n+1}(p)}{t_{n}(p)} \right)$ of different \(p\)-sequences are different, say $\Phi_p$, and tends towards $2$ for $p$ tending towards infinity. That is,
\be \label{eq-p-ratio-limit}
\lim_{n \rightarrow \infty} \frac{t_{n+1}(p)}{t_{n}(p)} = \Phi_p \overset{p \rightarrow \infty} \longrightarrow 2. \ee

Similarly, if we define another \(p\)-sequence as
\be \tilde{t}_{n+p}(p) := \left(t^m_{n+p-1}(p) + t^m_{n+p-2}(p) + \cdots + t^m_{n}(p) \right)^{\frac{1}{m}}, \ee 
then
\be \label{eq-p-ratio-limit1}
\lim_{n \rightarrow \infty} \left(\frac{\tilde{t}_{n+1}(p)}{\tilde{t}_{n}(p)} \right)^{\frac{1}{m}} = \left(\Phi_p \right)^{\frac{1}{m}}  \overset{m \rightarrow \infty} \longrightarrow 1. \ee

In particular, when $m=2$ and $p \rightarrow \infty$,we have
\be \label{eq-p-ratio-limit2}
\lim_{n \rightarrow \infty} \left(\frac{\tilde{t}_{n+1}(p)}{\tilde{t}_{n}(p)} \right)^{\frac{1}{2}} = \left(\Phi_p \right)^{\frac{1}{2}}  \overset{p \rightarrow \infty} \longrightarrow \sqrt{2}. \ee
Note also that $\sqrt{2}$ is the ratio of the diagonal of a square to its side.
Furthermore, if one considers a square of side length $a$, and its diagonal $\sqrt{2}a$ as the side of another square, and so on, then
\bea 
&& \frac{\textit{Side of square}}{\textit{Side of preceding square}} = \sqrt{2} = \frac{\textit{Diagonal of square}}{\textit{Side of square}}, \label{eq-sqdiag-ratio} \\
&& \frac{\textit{Area of square}}{\textit{Area of preceding square}} = 2. \label{eq-area-ratio}
\eea
Similarly, if one considers a cube of side length $a$, and its diagonal $\sqrt{3}a$ as the side of another cube, and so on, then
\bea 
&& \frac{\textit{Side of cube}}{\textit{Side of preceding cube}} = \sqrt{3} = \frac{\textit{Diagonal of cube}}{\textit{Side of cube}}, \label{eq-cubediag-ratio} \\
&& \frac{\textit{Volume of cube}}{\textit{Volume of preceding cube}} = (\sqrt{3})^3. \label{eq-vol-ratio}
\eea

\section{\(p\)-golden ratio}

\begin{figure}%
\centering
\includegraphics[width = 3in]{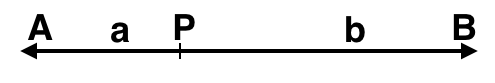}%
\caption{Division of a line into \(2\) segments.}
\label{fig-gratio2}
\end{figure}

The golden ratio \cite{livio2002} arises when we consider division of a line segment $AB$ with a point $P$ such that  $\frac{BP}{AP} = \frac{AB}{BP}$, where $BP > AP$
(see Fig. \ref{fig-gratio2}). Given $AP = a$ and $BP = b$ are two positive numbers, the above problem translates as
\be \label{eq-2ratio} \frac{b}{a} = \frac{a+b}{b}. \ee 
Taking $\frac{b}{a} = x$, the above equation can be rewritten as $x = 1 + \frac{1}{x}$. This reduces to the characteristic equation 
\be \label{eq-2ce} X(x) = x^2 - x -1=0, \ee 
whose positive solution is 
\be \Phi = \frac{\sqrt{5}+1}{2} = 1.618 \cdots. \ee 
The conjugate of $\Phi$ is its reciprocal, given by
\be \varphi = \frac{\sqrt{5}-1}{2} = \frac{1}{\Phi} = 0.618 \cdots. \ee
The golden ratio allegedly appears everywhere: in geometry, math, science, art, architecture, nature, human body, music, painting. However, many hold skeptical views about this \cite{devlin, falbo1, gould1981, markowsky1992, spira2012}.\\

\begin{figure}[htb]%
\centering
\includegraphics[width = 3in]{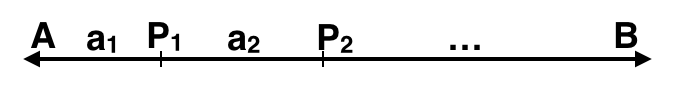}%
\caption{Division of a line into \(p\) segments.}
\label{fig-gratiop}
\end{figure}

We wish to generalize the above case. That is, we consider division of a line segment $AB$ into \(p > 2\) segments (see Fig. \ref{fig-gratiop}) such that 
$AP_1(= a_1) < P_1P_2 (= a_2) < \cdots < BP_{p-1} (=a_p)$ are \(p \) positive real numbers. 
We now demand that
\bea 
&& \frac{P_1P_2}{AP_1} = \frac{P_2P_3}{P_1P_2} = \cdots = \frac{AB}{P_{p-1}B} \label{eq-p-ratio1} \\
&& \Leftrightarrow \frac{a_2}{a_1} = \frac{a_3}{a_2} = \cdots = \frac{\sum_{k=1}^p a_k}{a_p}.  \nonumber 
\eea
If a unique positive value, say $\Phi_p$, exists for the above ratio, we call it the \(p\)-golden ratio.
%
From Eq. (\ref{eq-p-ratio1}) follows naturally the $p$-degree algebraic equation whose {\it positive solution} gives the value of $\Phi_p$ \cite{asutosh-gr}:  
\be \label{eq-godeqn}
X_p(x) \equiv x^p - \sum_{k=1}^{p-1} x^k - 1 = 0. \ee
Note that $X_p(0) = -1$ for all \(p\) and $X_p(1) = -(p-1)$.
Eq. (\ref{eq-godeqn}) is the {\it characteristic equation} for $\Phi_p$.
Interestingly, {\it $\Phi_p$ coincides with the limiting ratio value} of the \(p\)-sequence \cite{asutosh-gr}.

\section{Relations between \(e\), \(\pi\), $\Phi_p$ and $\sqrt{2}$}
In this section, firstly we present the relations between \(e\), \(\pi\), the \(2\)-golden ratio, and \(\sqrt{2}\).  
The Fibonacci sequence $\{0,1,1,2,3,5,8,13,21,34,55,\cdots \}$ is a \(2\)-sequence because it is generated by the sum of two previous terms. 
The positive real solution of the characteristic equation $x^2 - x -1 =0$ yields the \(2\)-golden ratio,
\be \label{eq-2-ratio-limit}
\lim_{n \rightarrow \infty} \frac{t_{n+1}(2)}{t_{n}(2)} = \Phi_2 = \frac{\sqrt{5} + 1}{2}. \ee

\noindent Following relations hold between \(e\), \(\pi\), $\Phi_2 \equiv \Phi$, and $\sqrt{2}$.
\bea
&& \Phi = 2\cos 36^{\circ} = e^{i \pi/5} + e^{-i \pi/5}, \label{eq-rel2} \\
&& \Phi (\Phi - 1) = e^{i2\pi} = - e^{i\pi} = -i e^{i \pi/2}, \label{eq-rel3} \\
&& \Phi (e^{i\pi} + \Phi) = 1, \label{eq-rel5} \\
&& \Phi \approx \sqrt{2} \left(1+ \frac{\sqrt{2}}{10} \right), \\
&& \pi \approx \sqrt{2} + \sqrt{3} - \frac{1}{100 \Phi}, \label{eq-b} \\
&& \frac{1}{e} + \frac{1}{\Phi} = 0.9859 \cdots \approx 1, \\
&& e \Phi \approx e + \Phi + \frac{1}{10 \Phi}, \\
&& e \approx \Phi^2 +0.1 = 2.7180339887, \label{eq-rel1} \\ 
&& \pi^2 \approx e^2 + \Phi^2 - \frac{\pi e \Phi}{100}, \\
&& \frac{\pi + e}{(\sqrt{3})^3 + \frac{1}{\Phi}} = 1.007858 \cdots, \label{eq-a} \\
&& \frac{\pi} {e} - \frac{\Phi} {\sqrt{2}} =  0.0116 \cdots, \label{eq-rel7} \\
&& \frac{\pi/ e}{\Phi / \sqrt{2}} = 1.0101427 \cdots, \label{eq-rel8} \\
&& \sqrt{2} \approx \frac{\sqrt{2} + \sqrt{3}}{\sqrt{5}} = \frac{\sqrt{2} + \sqrt{3}}{5} \left(\Phi + \frac{1}{\Phi}\right) \approx \frac{e \Phi} {\pi} \approx \frac{7}{5}. \label{eq-rel9}
\eea

A transcendental expression cannot be equated with an algebraic expression.
Eq. (\ref{eq-rel1}), however, give the polynomial approximation of $e$ in terms of $\Phi$. \\

\noindent 
Because $\Phi_p$ is a solution of Eq. (\ref{eq-godeqn}), we have 
\bea 
\Phi_p^p &=& \Phi_p^{p-1} + \Phi_p^{p-2} + \cdots + \Phi_p + 1 = \sum_{k=0}^{p-1} \Phi_p^k, 
\label{eq-basic} \\
\Phi_p^{p+1} &=& \Phi_p^{p} + \Phi_p^{p-1} + \cdots + \Phi_p^2 + \Phi_p, \nonumber \\
&=& 2\Phi_p^{p} - 1. \label{eq-basic0}
\eea
Using these equations, we have the following relations between \(e\), \(\pi\) and $\Phi_p$ ($p \ge 3$).
\bea 
&&  e^{i\pi} + \Phi_p^p - \sum_{k=1}^{p-1} \Phi_p^k = 0, \label{eq-rel-p1} 
\eea

Eq. (\ref{eq-rel-p1}) has been obtained using the Euler's identity.

\section{Conclusion}
In summary, inspired by Euler's identity, we have provided several relations between the base of natural logarithms ($e$), the ratio of the circumference of a circle to its diameter ($\pi$), the golden ratios ($\Phi_p$) of the additive \(p\)-sequences, and the ratio of the diagonal of a square to its side ($\sqrt{2}$).

\vspace{0.5cm}
\noindent {\Large \bf Acknowledgements}\\
AK would like to thank all the readers for their useful comments which helped to improve the manuscript, and Bruce E. Camber for motivating to add $\sqrt{2}$ to the list of $\{e, \pi, \Phi_p\}$.


\end{document}